\documentclass[10pt]{article}

\usepackage[reqno]{amsmath}
\usepackage{amssymb, amsthm, enumerate, url}

\newcommand\de{\delta}

\newcommand\la{\lambda}

\renewcommand\th{\theta} 

\newcommand\re{{\mathbb R}}

\newcommand\opk[1]{\mathop{\mathrm{#1}}\nolimits}

\newcommand\comp[1]{{\mkern2mu\overline{\mkern-2mu#1}}}
\newcommand\diff{\mathbin{\mkern-1.5mu\setminus\mkern-1.5mu}}

\newcommand\inv{^{-1}}
\newcommand\sbs{\subseteq}

\newcommand\seq[3]{#1_{#2},\ldots,#1_{#3}}

\newcommand\pmat[1]{\begin{pmatrix} #1 \end{pmatrix}}

\newcommand\one{{\bf1}}

\newcommand\col{\opk{col}}

\newtheoremstyle{plainsl}%
     {\topsep}
     {\topsep}
     {\slshape} 
     {}
     {\normalfont\bfseries}
     {.}
     { }
     {}

\swapnumbers

\theoremstyle{plainsl}
\newtheorem{theorem}{Theorem}[section]
\newtheorem{lemma}[theorem]{Lemma}
\newtheorem{corollary}[theorem]{Corollary}

\renewcommand\proof{\noindent\textsl{Proof. }}
\newcommand\sqr[2]{{\vbox{\hrule height.#2pt
     \hbox{\vrule width.#2pt height#1pt \kern#1pt
         \vrule width.#2pt}\hrule height.#2pt}}}
\renewcommand\qed{
     \ifmmode\eqno\sqr53
     \else\nobreak\ \hfill\sqr53\medbreak\fi}

\newcommand\lref[1]{Lemma~\ref{lem:#1}}
\newcommand\tref[1]{Theorem~\ref{thm:#1}}
\newcommand\cref[1]{Corollary~\ref{cor:#1}}
\newcommand\sref[1]{Section~\ref{sec:#1}}

\newcommand\bh[2]{#1^{\{#2\}}}
\newcommand\bhsx{\bh{X}{2}}
\newcommand\cprod{\mathbin{\scriptscriptstyle\square}}
\DeclareMathOperator{\dist}{dist}

\newcommand\dip[2]{\langle#1,#2\rangle}
\newcommand{\trace}{\mathop{\rm Tr}\nolimits}

\newcommand{\bra}[1]{#1^*}
\newcommand{\ket}[1]{#1}

\newcommand{\identity}{{\mathbb{I}}}
\newcommand{\F}{{\mathbb{F}}}
\newcommand{\J}{{\mathbb{J}}}
\newcommand{\C}{{\mathbb{C}}}

\newcommand{\bea}{\begin{eqnarray}}
\newcommand{\eea}{\end{eqnarray}}
\newcommand{\beas}{\begin{eqnarray*}}
\newcommand{\eeas}{\end{eqnarray*}}

\begin{document}

\begin{center}
{\Large Symmetric Squares of Graphs}

\bigskip
{\large
Koenraad Audenaert$\mbox{}^{1,4}$, Chris Godsil$\mbox{}^2$, Gordon
Royle$\mbox{}^3$ and Terry Rudolph$\mbox{}^{1,4}$}

\medskip
${}^1$ Department of Physics, Blackett Laboratory, Imperial College
London,\\ Prince Consort Road, London SW7 2BW, UK \\
${}^2$ Department of Combinatorics and Optimisation, University of
Waterloo,\\ Waterloo, Ontario N2L 3G1, Canada \\
${}^3$ Department of Computer Science \& Software Engineering,\\
University of Western
Australia, Crawley, WA 6009, Australia\\
${}^4$ Institute for Mathematical Sciences, Imperial College London,\\
53 Exhibition Road, London SW7 2BW, UK
\end{center}
\begin{abstract}
We consider \textsl{symmetric powers} of a graph. In particular, we
show that the spectra of the symmetric square of strongly regular
graphs with the same parameters are equal. We also provide some
bounds on the spectra of the symmetric squares of more general
graphs. The connection with generic exchange Hamiltonians in quantum
mechanics is discussed in an appendix.
\end{abstract}

\section{Introduction}
The \textsl{symmetric $k$-th power} $\bh{X}{k}$ of a graph $X$ is
constructed as follows:  its vertices are the $k$-subsets of $V(X)$,
and two $k$-subsets are adjacent if and only if their symmetric
difference is an edge.  As an example, and a test case, the
symmetric square of the complete graph $K_n$ is its line graph.
(Useful procedures for constructing symmetric squares of arbitrary
graphs will be given in Theorem \ref{symmsq} and Lemma
\ref{symmsq2}). Symmetric powers were introduced in \cite{Rud01}.

The symmetric powers are related to a class of random walks, where
one starts with $k$ particles occupying $k$ distinct vertices of
$X$, and, at each step of the walk, a single particle moves to an
unoccupied adjacent site.  More formally, we can generalise
the concept of a walk on a graph to a $k$-walk, which is an
alternating sequence of $k$-subsets of vertices $V_i$ and arcs
$e_i$, $(V_0,e_1,V_1,e_2\ldots,e_n,V_n)$, such that the symmetric
difference of $V_{i-1}$ and $V_i$ is the arc $e_i$. It is readily
seen that a $k$-walk on $X$ corresponds to an ordinary 1-walk on
$\bh{X}{k}$.

Our motivation for studying symmetric powers arises from its relevance
for physically realisable systems and for the graph isomorphism problem.
A brief outline of the connection between symmetric powers and
exchange Hamiltonians in quantum mechanics is given in the appendix.

The relevance to the graph isomorphism problem arises because
invariants of the symmetric powers of $X$ are invariants of $X$.
There are examples of cospectral graphs $X$ and $Y$ such that
$\bh{X}2$ and $\bh{Y}2$ are not cospectral. In fact we have
verified computationally that graphs on at most 10 vertices
determined by the spectra of their symmetric squares. On the other
hand, the main result of this paper is a proof that if $X$ and $Y$
are cospectral strongly-regular graphs then $\bh{X}2$ and $\bh{Y}2$
are cospectral.  There is also a family of five regular graphs on 24
vertices whose symmetric squares are cospectral.  Nevertheless, in
each of those cases, and, in fact, for all graphs we have examined
(including strongly regular graphs on up to 36
vertices), the spectrum of the symmetric \textsl{cubes} determine
the original graphs.  (The computations on the the strongly regular
graphs on 35 and 36 vertices were performed by Dumas, Pernet and
Saunders; more details are given in Section~\ref{comp-res}.)

If it were true for some fixed $k$ that any two graphs $X$ and $Y$ are
isomorphic if and
only if their $k$-th symmetric powers are cospectral, then we would
have a polynomial-time algorithm for solving the graph isomorphism
problem. For a
pessimist this suggests that, for each fixed $k$, there should be
infinitely many pairs of non-isomorphic graphs $X$ and $Y$ such that
$\bh{X}k$ and $\bh{Y}k$ are cospectral.


In the last section of the paper we will consider bounds, from an
algebraic perspective, on the spectra of the symmetric squares of
arbitrary graphs.

While the focus of this paper is on the spectra of the symmetric
squares, it should be noted that multivalued graph invariants based
on generic (analytic) matrix valued functions $f(\bh{A}k)$ can also
be considered \cite{Rud01}, where $\bh{A}k$ is the adjacency matrix
of $\bh{X}k$. In \cite{SJC} this approach was followed, and
numerical computations showed that the values of $\exp(i\bh{A}{2})$
sufficed to distinguish all strongly regular graphs up to around 30
vertices.


\section{Preliminaries}

If $A$ is square matrix, then let $\phi(A,t)$ denote the
characteristic polynomial $\det(t\identity-A)$ of $A$.  If $A$ is the
adjacency matrix of $X$, we will also write $\phi(X,t)$.  If $x$ and
$y$ are vertices of $X$, we write $x\sim y$ to denote that $x$ is
adjacent to $v$.

A graph is \textsl{strongly regular} with parameters $(v,k;a,c)$ if it
is not complete or empty,
has $v$ vertices, and the number of common neighbours of two vertices
$x$ and $y$ is $k$,
$a$ or $c$ according as $x$ and $y$ are
equal, adjacent, or distinct and not adjacent.  Thus if $X$ is strongly
regular, the neighbourhood
of each vertex in $X$ is regular and the neighbourhood of each vertex
in the complement of
$X$ is regular.
The line graph of the complete graph $K_n$ is strongly regular if
$n\ge3$.

The main tool in this paper will be walk-generating functions.  If $A$
is the adjacency matrix of
the graph $X$, then the walk-generating function $W(X,t)$ is the formal
power
series
\[
\sum_{r\ge0} A^r t^r.
\]
We view this either as a power series with coefficients from the ring
of matrices, or
as a matrix whose entries are power series over $\re$.  Its $ij$-entry
$W_{i,j}(X,x)$
is the generating function for the walks in $X$ that start at the
vertex $i$ and finish at $j$.

If $D\sbs V(G)$, then $W_{D,D}(X,t)$ denotes the submatrix of
$W(X,t)$ with rows and columns indexed by the vertices in $D$.  The
following identities are proved in Chapter 4 of \cite{God}.

\begin{theorem}
\label{thm:deld}
If $D$ is a subset of $d$ vertices of $X$, then
\[
t^{-d}\det(W_{D,D}(X,t\inv)) =\frac{\phi(X\diff D,t)}{\phi(X,t)}.\qed
\]
\end{theorem}

\begin{corollary}
\label{cor:delvx}
If $i\in V(X)$, then
\[
t\inv W_{i,i}(X,t\inv) =\frac{\phi(X\diff i,t)}{\phi(X,t)}.\qed
\]
\end{corollary}

\begin{corollary}
\label{cor:del2}
If $i$ and $j$ are distinct vertices of $X$,
\[
t\inv W_{i,j}(X,t\inv)
     =\frac{(\phi(X\diff i,t)\phi(X\diff j,t)-\phi(X,t)\phi(X\diff
ij,t))^{1/2}}{\phi(X,t)}.\qed
\]
\end{corollary}

The presence of the square root in the previous identity is surprising.
  Note though
that it causes no ambiguity, since we know that the coefficients of
$W_{i,j}(G,t)$
are non-negative.

We apply these identities to obtain information about strongly
regular graphs.  If $X$ is strongly regular with parameters
$(v,k;a,c)$ and adjacency matrix $A$ then
\[
A^2-(a-c)A-(k-c)\identity =c\J.
\]
(This is essentially the definition of ``strongly regular''
expressed in linear algebra.) Since $X$ is regular $A$ and $\J$
commute, whence we see that for each non-negative integer $k$, the
power $A^k$ is a linear combination of $\identity$, $\J$ and $A$. Thus
the
generating function $W_{i,j}(X,t)$ depends only on whether the
vertices $i$ and $j$ are equal, adjacent, or distinct and not
adjacent.  Using the corollaries above, this leads to the following:

\begin{theorem}
Let $X$ be a strongly regular graph.
Then $\phi(X\diff i,t)$ is independent of $i$ and, if $i\ne j$, then
$\phi(X\diff ij,t)$ only depends on whether $i$ and $j$ are adjacent or
not. \qed
\end{theorem}

\begin{theorem}
Let $X$ be a strongly regular graph and let $D_1$ and $D_2$ be induced
subgraphs of
$V(X)$. If $D_1$ and $D_2$ are cospectral with cospectral complements,
then
$X\diff D_1$ and $X\diff D_2$ are cospectral with cospectral
complements.
\end{theorem}

\proof
Suppose $D\sbs V(X)$.  Then $W_{D,D}(X,t)$ is the submatrix of $W(X,t)$
with rows and columns indexed by the vertices in $D$.  Since $X$ is
strongly
regular, we have
\[
W_{D,D}(X,t) =\alpha\identity+\beta\J + \gamma A(D)
\]
where $\alpha$, $\beta$ and $\gamma$ are generating functions
and $A(D)$ is the adjacency matrix of the subgraph induced by $D$. So
\begin{align*}
\det(\alpha\identity+\beta\J+\gamma A(D))
&=\det((\alpha\identity+\gamma A(D))(\identity+(\alpha\identity+\gamma
A(D))\inv \beta\J))\\
&=\det(\alpha\identity+\gamma
A(D))\det(\identity+(\alpha\identity+\gamma A(D))\inv \beta\J)
\end{align*}

Recall that if the matrix products $BC$ and $CB$ are defined then
\[
\det(\identity+BC) =\det(\identity+CB).
\]
Since
\[
=\det(\identity+(a\identity+cA(D))\inv b\one\one^T)
\J=\one\one^T
\]
it follows that
\[
\det(\identity+(\alpha\identity+\gamma A(D))\inv \beta\J)
=1+\beta\one^T(\alpha\identity+\gamma A(D))\inv\one.
\]
We are working effectively over the field of real rational functions in
$t$,
therefore
\[
\det(\alpha\identity+\gamma A(D)) =
\gamma^{|D|}\det\left(\frac{\alpha}{\gamma}\identity+A(D)\right)
\]
and
\[
\one(\alpha\identity+\gamma A(D))\inv\one^T
     =\alpha\inv\sum_{r\ge0} \left(\frac{\gamma}{\alpha}\right)^r
\dip{\one}{A^r\one}.
\]
We conclude that $\det W_{D,D}(X,t)$ is determined by
\[
\alpha,\ \beta,\ \gamma,\ \phi(D,t)
\]
and the series
\[
\sum_{r\ge0} \dip{\one}{A(D)^r\one} \, t^r
\]
which is the generating function for all walks in $D$.  By Exercise~10
in
Chapter 4 of \cite{God}, this generating function is determined by
the characteristic polynomial of $D$ and its complement.

Consequently we have shown that if $D_1$ and $D_2$ are induced
subgraphs of $X$,
cospectral with cospectral complements, then $X\diff D_1$ and $X\diff
D_2$ are cospectral.
Applying this to the complement of $X$, which is also strongly regular,
we deduce that
the complements of $X\diff D_1$ and $X\diff D_2$ are cospectral.
\qed

If $S_1$ and $S_2$ are independent sets of the same size in the
strongly regular graph $X$,
the previous theorem implies that $X\diff S_1$ and $X\diff S_2$ are
cospectral.  Even this
special case of the theorem appears to be new.

\section{Equitable Partitions}
\label{sec:eqpars}

We will also be working with equitable partitions of graphs.  A
partition $\pi$ of the
vertices of $X$ is \textsl{equitable} if for each pair of cells $C_i$
and $C_j$ of $\pi$
there is constant $b_{i,j}$ such that each vertex in $C_i$ has exactly
$b_{i,j}$
neighbours in $C_j$.  The {\textsl{quotient graph} $X/\pi$ has the
cells of $\pi$ as its
vertices, with $b_{i,j}$ directed edges from $C_i$ to $C_j$.  If $G$ is
a group of
automorphisms
of $X$, then the orbits of $G$ form an equitable partition.  If $X$ is
strongly regular and
$u\in V(X)$, the partition with three cells consisting of $\{u\}$, the
neighbours of $u$,
and the vertices at distance two from $u$ is equitable.

If $\pi$ is a partition, the \textsl{characteristic matrix} of $\pi$
is the matrix with the characteristic vectors of the cells of $\pi$
as its columns.  (Thus it is a 01-matrix and each row-sum is equal
to 1.)  If $\pi$ is an equitable partition of $X$ with
characteristic matrix $R$ and $B:=A(X/\pi)$, then
\[
AR =RB.
\]
There is a matrix $B$ such that $AR=RB$ if and only if $\col(R)$ is
$A$-invariant,
and this in turn holds if and only if $\pi$ is equitable.  If $z$ is an
eigenvector for
$B$ with eigenvalue $\la$, then $Rz$ is an eigenvector for $A$ with
eigenvalue $\la$.  This
shows that each eigenvalue of $B$ is an eigenvalue of $A$.

As a particularly relevant example, the symmetric
square $\bhsx$ has two sorts of vertices: the pairs $uv$ where $u\sim
v$ and the pairs
$uv$ where $u\not\sim v$.  If $X$ is strongly regular with parameters
$(v,k;a,c)$,
then this partition is equitable with quotient matrix
\[
B = \pmat{2a & 2k-2a-2\\ 2c & 2k-2c}.
\]
If $\de := a-c$, then the eigenvalues of this matrix are
\[
k+\de\pm\sqrt{(k-\de)^2-4c},
\]
and these are eigenvalues of the symmetric square.  The eigenvector
$z$ of $B$ corresponding to the positive eigenvalue if positive, and
therefore $Rz$ is a positive eigenvector of $A$. This implies that
the positive eigenvalue is the spectral radius of the symmetric
square.

We have the following relation between walks in $X$ and $X/\pi$ when
$\pi$ is
equitable.

\begin{lemma}
\label{lem:qwalks}
Let $X$ be a graph with adjacency matrix $A$.  If $\pi$ is an equitable
partition of $X$
and $B:=A(X/\pi)$, then the $r,s$-entry of $B^k$ is equal to the number
of walks of
length in $X$ that start on a given vertex in cell $C_r$ and finish on
a vertex on $C_s$.
\end{lemma}

\proof
Assume $v=|V(X)|$.
Let $\pi$ be an equitable partition of $X$ with $r$ cells and let $R$
be the characteristic
matrix of $\pi$.  Then $AR=RB$ and, more generally,
\[
A^kR =RB^k,\quad k\ge0.
\]
Let $\seq e1v$ denote the standard basis of $\re^v$ and let $\seq f1r$
denote the standard
basis of $\re^r$.  Let $u$ and $v$ be vertices of $X$ that form
singleton cells of $\pi$,
and suppose $\{v\}$ is the $j$-th cell of $\pi$.  If $u\in V(X)$ then
\[
\dip{e_U}{A^kRf_j}
\]
is the number of walks of length $k$ in $X$ that start at $u$ and
finish on a vertex in
the $j$-th cell of $\pi$.  On the other hand, if vertex $u$ is in the
$i$-th cell of $\pi$,
then $Re_u=f_i$ and
\[
\dip{e_u}{RB^kf_j} =\dip{f_i}{B^kf_j}.\qed
\]

\section{Constructing the Symmetric Square}

The main result of this paper depends on the observation that we can
construct
the symmetric square of $X$ in two stages.

We begin with the Cartesian product of $X$ with itself, which has
adjacency matrix
\[
A\otimes \identity +\identity\otimes A.
\]
The vertex set of the Cartesian product $X\cprod Y$ of $X$ and $Y$ is
$V(X)\times V(Y)$,
and $(x,y)\sim (x',y')$ if either $x=x'$ and $y\sim y'$ , or $x\sim x'$
and $y=y'$.
We also have
\[
\dist_{X\cprod Y}((x,y),(x',y')) =\dist_X(x,x')+\dist_Y(y,y').
\]
We denote $X\cprod X$ by $X^{\cprod2}$.
The subgraph of $X^{\cprod2}$ induced by the vertices
\[
\{(i,i): i\in V(X)\}
\]
is called the \textsl{diagonal}.

The map
\[
\tau:(i,j) \mapsto (j,i)
\]
is an automorphism of $X^{\cprod2}$.  It fixes each vertex in the
diagonal and partitions the remaining vertices into pairs.  We will
call it the \textsl{flip} automorphism of $X^{\cprod2}$.

\begin{theorem}\label{symmsq}
Let $X$ be a graph, let $D$ denote the diagonal of $X^{\cprod2}$ and
let $\pi$ be the
partition of $(X^{\cprod2})\diff D$ formed by the non-trivial orbits of
the flip.  Then
$\bhsx$ is isomorphic to $((X^{\cprod2})\diff D)/\pi$.\qed
\end{theorem}

We make some comments on the quotienting involved.  Suppose $i$ and
$j$ are distinct vertices in $X$.  Then $(i,j)\not\sim(j,i)$, and
therefore each orbit of the flip of size two is an independent set.
If $i\ne\ell$ and $(i,j)\sim(i,\ell)$, then $(i,j)\not\sim(\ell,i)$.
Hence two orbits of the flip are either not joined by any edges, or
else each vertex in one orbit has exactly one orbit in the second.
It follows from this that $((X^{\cprod2})\diff D)/\pi$ has no loops
and no multiple edges---it is a simple graph.

Our aim now is to show that if $X$ and $Y$ are strongly regular graphs
with the same
parameters, then the graphs obtained by deleting the diagonal from
$X^{\cprod2}$ and
$Y\cprod Y$ are cospectral (with cospectral complements).  We will then
show that
the quotients modulo the flip are cospectral.

\section{Deleting the Diagonal}

If $\th$ is an eigenvalue of $A$, let $E_\th$ denote the orthogonal
projection onto
the eigenspace belonging to $\th$.  Then if $r\ge0$, we have the
spectral decomposition:
\[
A^r =\sum_\th \th^r E_\th.
\]
from which we have
\[
W(X,t) =\sum_\th (1-t\th)\inv E_\th.
\]
Since
\[
A\otimes \identity+\identity\otimes A
=\sum_{\th,\tau}(\th+\tau)E_\th\otimes E_\tau,
\]
we see that
\[
W(X^{\cprod2},t) =\sum_{\th,\tau}(1-t(\th+\tau))\inv E_\th\otimes
E_\tau.
\]

If $M$ and $N$ are $m\times n$ matrices, their \textsl{Schur product}
(also called Hadamard product)
$M\circ N$ is the $m\times n$ matrix given by
\[
(M\circ N)_{i,j} =M_{i,j}N_{i,j}.
\]

\begin{theorem}
\label{thm:wgschur}
If $D$ denotes the diagonal of $X^{\cprod2}$ and $A(X)$ has the
spectral decomposition
$\sum_\th \th E_\th$, then
\[
W_{D,D}(X^{\cprod2},t) =\sum_{\th,\tau}(1-t(\th+\tau))\inv E_\th\circ
E_\tau.
\]
\end{theorem}

\proof
It is enough to note that
\[
(E_\th\otimes E_\tau)_{D,D} =E_\th\circ E_\tau.\qed
\]

The linear span of the principal idempotents of the adjacency matrix of
a strongly
regular graph is equal to the span of $\identity$, $A(X)$ and $\J$, and
is therefore closed under
the Schur product.  Hence $E_\th\circ E_\tau$ is a linear combination
of principal idempotents.
The coefficients in this linear expansion are known as the
\textsl{Krein parameters} of
the strongly regular graph, and are determined by the parameters of the
graph.
Therefore the eigenvalues of $W_{D,D}(X^{\cprod2},t)$ are determined
by the parameters of $X$, and so $\det(W_{D,D}(X^{\cprod2},t) )$ is
determined
by the parameters of $X$.

\begin{lemma}
\label{lem:cut}
If $X$ is a strongly regular graph and $D$ is the diagonal of
$X^{\cprod2}$, then the
spectrum of $X^{\cprod2}\diff D$ is determined by the spectrum of
$X$.\qed
\end{lemma}

\section{Flipping Quotients}
\label{sec:flip}

We use $Y$ to denote the quotient of $X^{\cprod2}$ by the flip.
By \lref{qwalks} we have the following.

\begin{lemma}
If $Y$ denotes the quotient of $X^{\cprod2}$ by the flip and $D$
denotes both the
diagonal of $X^{\cprod2}$ and the image of $D$ in $Y$, then
\[
\frac{\phi(Y\diff D,t)}{\phi(Y,t)}
     =\frac{\phi(X^{\cprod2}\diff D,t)}{\phi(X^{\cprod2},t)}.
\]
\end{lemma}

We now show that, for any graph $X$, the spectrum of $Y$ is determined
by the spectrum
of $X$. Given the above lemma it follows immediately that if $X$ is
strongly regular,
then the spectrum of $\bhsx$ is determined by the spectrum of $X$.

Let $X_1$ and $X_2$ be two cospectral graphs on $v$ vertices with
adjacency matrices
$A_1$ and $A_2$.
Let $L$ be an orthogonal matrix such that
\[
L^TA_1L =A_2.
\]

Let $\F$ be the permutation matrix that represents the flip on
$\re^v\otimes\re^v$.
So $\F$ maps $x\otimes y$ to $y\otimes x$, for all $x$ and $y$ in
$\re^v$.
Let $R$ be the normalized characteristic matrix of the orbit partition
of the flip---$R$
is obtained from the characteristic matrix of the orbit partition by
normalizing each column.
We have
\[
R^TR=\identity,\quad RR^T=\frac{1}{2}(\identity+\F).
\]
Let $A_i^{\cprod2}$ denote the adjacency matrix of $X_i^{\cprod2}$.
Then there
are matrices $C_i$ such that
\[
A_i^{\cprod2} R=RC_i
\]
We prove that $C_1$ and $C_2$ are cospectral.

We have
\[
C_2 =R^TA_2^{\cprod2}R =R^T(L\otimes L)^TA_1^{\cprod2}(L\otimes L)R
\]
whence
\[
RC_2R^T =RR^T(L\otimes L)^TA_1^{\cprod2}(L\otimes L)RR^T.
\]
Because $L\otimes L$ and $\F$ commute, so do $L\otimes L$ and $RR^T$.
So
\[
RC_2R^T =(L\otimes L)^TRR^T A_1^{\cprod2}RR^T(L\otimes L)
     =(L\otimes L)^TR\,C_1\,R^T(L\otimes L)
\]
and hence
\[
C_2 =R^T(L\otimes L)^TR\,C_1\,R^T(L\otimes L)R.
\]
Since
\begin{align*}
R^T(L\otimes L)^TR R^T(L\otimes L)R &=R^T(L\otimes L)^T (L\otimes
L)RR^TR\\
     &=R^TRR^TR\\
     &=\identity,
\end{align*}
we conclude that $C_1$ and $C_2$ are similar matrices.

\medbreak
Note that it is possible to express the spectrum of $Y$ in terms of the
spectrum of $X$.
If $\pi$ is equitable and $B=A(X/\pi)$ and $\th$ is an eigenvalue of
$B$, then
\[
\dim(\ker(B-\th \identity)) = \dim\bigl(\col(R)\cap\ker(A-\th
\identity)\bigr).
\]
Suppose $\seq z1n$ is an orthonormal basis for $\re^n$ consisting of
eigenvectors
of $X$. Then the products $z_i\otimes z_j$ form an orthonormal basis
for $\re^{n^2}$
consisting of eigenvectors of $X^{\cprod2}$. If $i\ne j$ then the span
of $z_i$ and $z_j$ is equal to the span of
the symmetric and antisymmetric combinations
\[
(z_i\otimes z_j)+(z_j\otimes z_i),\quad (z_i\otimes z_j)-(z_j\otimes
z_i)
\]
These two vectors are orthogonal and the first is constant on the orbit
partition of the
flip, while the second sums to zero on each orbit.  If $z^TA=\th z$
then $z^TRB=\th z^TR$.
So if $\th$ has multiplicity $\ell$ as an eigenvalue of $X$, the vectors
\[
(z_i\otimes z_j)+(z_j\otimes z_i),\quad z_i\otimes z_i,
\]
where $z_i\in\ker(A-\th \identity)$, give rise to a subspace of
eigenvectors of $Y$
with eigenvalue $2\th$ and dimension $\binom{\ell+1}2$.  If $\th$ has
multiplicity
$\ell$ and $\tau$ has multiplicity $m$, then we obtain a subspace of
eigenvectors of
the quotient with dimension $\ell m$.  By adding up the dimensions of
these subspaces,
we find that the images of the given vectors provide a basis consisting
of eigenvectors
of $Y$.  It follows that the multiplicities of the eigenvalues of
$Y$ are determined by the eigenvalues of $X$ and their multiplicities.
(If $X$ has
exactly $r$ distinct eigenvalues, then $X^{\cprod2}$ has at most
$\binom{r+1}2$;
if $X^{\cprod2}$ has fewer eigenvalues, then the procedure just
described  will give the
multiplicities of the eigenvalues of $Y$, but does not lead to a simple
formula.)

\section{More Cospectral}

We have seen that if $X$ and $Y$ are strongly regular graphs with the
same parameters,
then their symmetric squares are cospectral.  Here we extend this.

\begin{lemma}
If $X$ and $Y$ are strongly regular graphs with the same parmeters,
then the
complements of their symmetric squares are cospectral.
\end{lemma}

\proof
> From Exercise 22 in Chapter 2 of \cite{God}, we have
\[
\phi(\comp{X},t+1) =(-1)^v \phi(X,t)(1-\one^T(tI+A)\inv\one).
\]
> From this it follows that cospectral graphs $X$ and $Y$ have
> cospectral complements
if and only if the generating function for all walks in $X$ is equal to
the corresponding
generating function for $Y$.

Assume $X$ is strongly regular, let $A$ denote the adjacency matrix of
$\bhsx$,
let $\pi$ be the partition of the vertices of $\bhsx$ by valency and
let the
characteristic matrix $R$ and quotient matrix $B$ be defined as in
\sref{eqpars}.
Then $AR=RB$ and so, for if $\ell\ge0$,
\[
A^\ell R=RB^\ell.
\]
Since the columns of $R$ sum to $\one$,
\[
\one^T A^\ell\one^T =\one^TA^\ell R\one =\one^T RB^\ell\one.
\]
We have
\[
\one^T R =(vk/2,\ v(v-1-k)/2)
\]
and therefore the entries of $R^TB^\ell$ are determined by $\ell$ and
the parameters of $X$.
Hence the generating function for all walk in $\bhsx$ is determined by
the parameters of
the strongly regular graph $X$, and the result follows.\qed

\section{Variations}

The \textsl{direct product} $X\times Y$ of graphs $X$ and $Y$ has
vertex set equal to
$V(X)\times V(Y)$, and $(u,v)\sim(x,y)$ if and only if $u\sim x$ and
$v\sim y$.
We have
\[
A(X\times Y) =A(X)\otimes A(Y).
\]
The flip map
\[
(x,y) \mapsto (y,x)
\]
is again an automorphism of $X\times X$ that fixes the diagonal.
We can obtain an analog of the symmetric product by deleting the
diagonal and then
quotienting over the flip.  A slightly modified version of the argument
in this paper shows
that if $X$ is strongly regular, then the spectrum of this analog is
determined by the spectrum
of $X$.  The key step is to verify the following analog of
\tref{wgschur}:
\[
W_{D,D}(X^{\otimes2},t) =\sum_{\th,\tau}(1-t\th\tau)\inv E_\th\circ
E_\tau.
\]

For a second analog, we turn to the graph obtained from the Cartesian
power $X^{\cprod k}$
by deleting the diagonal and the quotienting over the orbits of the
automorphism
that sends each $k$-tuple to its right cyclic shift.  Again our
argument shows that
if $X$ is strongly regular, the spectrum of this analog is determined
by $X$.
Thus there is more than one candidate for the ``symmetric cube" of a
graph, but the
spectrum of the one just described is a less useful graph invariant
than the spectrum
of the symmetric cube defined in Section 1.

\section{Symmetric Squares of General Graphs}

In this section we take a closer look at the purely algebraic
properties of the symmetric powers, and of the symmetric square in particular.
We start by giving a purely algebraic definition.

Let $P^{(k)}$ be the 0/1-matrix with
${\binom{v}{k}}$ rows, labelled by the $k$-tuples
$(i,j,\ldots, l)$ with $1\leq i<j<\ldots<l\leq v$, and
$v^{k}$ columns, labelled by the $k$-tuples $[i^{\prime },j^{\prime
},\ldots, l^{\prime }]$
with $1\leq i^{\prime },j^{\prime },\ldots,l^{\prime}\leq v$,
such that the elements $P_{(i,j,\ldots, l),[i^{\prime},j^{\prime
},\ldots, l^{\prime }]}^{(k)}$ are 1 iff $(i,j,\ldots, l)$
is a permutation of $[i^{\prime },j^{\prime },\ldots, l^{\prime }]$.
Then
\begin{lemma}\label{symmsq2}
The adjacency matrix $\bh{A}{k}(X)$ of $\bh{X}{k}$ is
$$
\bh{A}{k}(X) = \tfrac{1}{(k-1)!}P^{(k)}~\left( A(X)\otimes
\identity_v^{\otimes k-1}\right) ~P^{(k)*}
$$
\end{lemma}

We focus on the symmetric square, and more generally on the
properties of the linear map
\[
\Omega: G\mapsto \Omega(G) \equiv \bh{G}2 = P^{(2)}(G\otimes
\identity)P^{(2)*}.
\]
Henceforth, we will write $P$ instead of $P^{(2)}$.

Because $\Omega$ is the composition of the two completely positive maps
\cite{choi}
$A\mapsto A\otimes\identity$ and $A\mapsto BAB^*$,
$\Omega$ is completely positive itself. In particular, $\Omega$
preserves positive semi-definiteness.
One easily checks
\bea
P P^* &=& 2\identity_{\binom{d}{2}} \label{eq:pps} \\
P^* P &=& \sum_{i,j=1}^d (E_{ii}\otimes E_{jj} + E_{ij}\otimes E_{ji})
- 2\sum_{i=1}^v E_{ii}\otimes E_{ii}, \label{eq:psp}
\eea
where $\{E_{ij}\}$ is the standard matrix basis.

The spectrum of a general Hermitian matrix and the spectrum
of its symmetric square have the same average value.
When $G$ is an adjacency matrix this obviously has no import, because
adjacency matrices are traceless.  However, in certain quantum mechanical contexts the map $\Omega$ is applied to Hamiltonians which are not traceless.

\begin{theorem}
For $G$ a $v\times v$ Hermitian matrix,
$$
\trace[G]/v = \trace[\bh{G}2]/\binom{v}{2}.
$$
\end{theorem}
\textit{Proof.}
The partial trace of $P^*P$ over the second tensor factor,
defined as $\trace[(X\otimes\identity)A]=\trace[X\,\trace_2[A]]$, yields
\beas
\trace_2[P^*P] &=& \sum_{i,j=1}^v (E_{ii} \trace[E_{jj}] + E_{ij}
\trace[E_{ji}])
- 2\sum_{i=1}^v E_{ii} \trace[E_{ii}] \\
&=& \sum_{i,j=1}^v (E_{ii} + E_{ij} \delta_{ij}) - 2\sum_{i=1}^v E_{ii}
\\
&=& (v-1)\identity_v.
\eeas
Therefore,
\beas
\trace[\Omega(G)] &=& \trace[P^*P(G\otimes\identity)] \\
&=& \trace[G\,\trace_2[P^*P]] \\
&=& (v-1)\trace[G].
\eeas
Dividing by $v(v-1)$ yields the statement of the Theorem.
\qed

\subsection{Comparison between the spectrum of a matrix and the
spectrum of its symmetric square}
For a Hermitian matrix $A$, we denote by $\lambda_k^\downarrow(A)$ its
$k$-th largest eigenvalue, counting multiplicities.
Likewise, $\lambda_k^\uparrow(A)$ is its $k$-th smallest eigenvalue.

We prove the following:
\begin{theorem}\label{eigtheorem}
For any non-negative positive semi-definite $v\times v$ matrix $G$, the following
relation holds, for $1\le m\le v$:
\[
\lambda_m^\downarrow(G) \le\lambda_m^\downarrow(\Omega(G)).
\]
\end{theorem}

\textit{Proof.}
We refer to \cite{hj} or \cite{bhatia} for the basic matrix analytical
concepts and theorems.

say it at all :-)
Focusing on a particular value of $m$, $1\le m\le v$, we need to show
\[
\lambda_m^\downarrow(G) \le \lambda_m^\downarrow(P(G\otimes
\identity)P^*),
\]
for all $G\ge 0$, or, equivalently,
\begin{equation}\label{eq:1}
\lambda_m^\downarrow(P(G\otimes \identity)P^*) \ge 1,
\end{equation}
for all $G\ge 0$ with $\lambda_m^\downarrow(G)=1$.

First note that one needs to prove this only for $G$ a partial
isometry of rank $m$. Indeed, for every $G\ge 0$ with
$\lambda_m^\downarrow(G)=1$,
there exists a partial isometry $B$ of rank $m$
such that $G\ge B$. As noted above, $\Omega$ is a completely positive
map, hence
$\Omega(G)\ge \Omega(B)$. By Weyl monotonicity we then have
$\lambda_m^\downarrow(\Omega(G))\ge \lambda_m^\downarrow(\Omega(B))$.
Thus
(\ref{eq:1}) follows for $G$ if it holds for $B$.

Let us write $B$ as $B=Q^* Q$, with $Q\in M_{m,v}(\mathbb{C})$ and
$Q Q^*=\identity_m$. Let $\ket{q_j}$ be the $j$-th column of $Q$. Thus
the
$\ket{q_j}$ are $v$ $k$-dimensional vectors and
$$
\sum_{j=1}^v \ket{q_j}\bra{q_j} =  \identity_m.
$$

The matrix $ P(Q^* Q\otimes \identity)P^* $ has the same non-zero
eigenvalues as
\[
(Q\otimes \identity)P^* P(Q^*\otimes \identity).
\]
Using the explicit form (\ref{eq:psp}), a short calculation shows that
$$
\lambda_m^\downarrow(\Omega(Q^* Q)) =
    \lambda_m^\downarrow(\identity +  A) = 1+\lambda_m^\downarrow(A),
$$
where $A$ is a $v\times v$ block matrix with blocks $A_{i,j}$ of
size $m\times m$ given by
$$
A_{i,j} = (1-2\delta_{ij}) \ket{q_j}\bra{q_i}.
$$
We have
$$
\sum_{i=1}^v A_{i,i} = -\identity_m.
$$

We have to show that $\lambda_m^\downarrow(A)\ge 0$. To that
purpose, consider the principal submatrix $A'$ of $A$ consisting of
the $2\times 2$ upper left blocks:
$$
A' = \left(\begin{array}{cc}
A_{11} & A_{12} \\
A_{21} & A_{22}
\end{array}\right)
= \left(\begin{array}{rr} -\ket{q_1}\bra{q_1} &  \ket{q_2}\bra{q_1}
\\ [2mm]
  \ket{q_1}\bra{q_2} & -\ket{q_2}\bra{q_2}
\end{array}\right).
$$
If we can prove that $\lambda_m^\downarrow(A')\ge 0$, this implies
$\lambda_m^\downarrow(A) \ge 0$ via eigenvalue interlacing.

When $m=1$, the $\ket{q_i}$ are scalars, and direct calculation
shows that $\lambda_1^\downarrow(A') = 0$.

For $m>1$, consider a (non-orthogonal) basis of $\mathbb{C}^d$ in
which $\ket{q_1}$ and $\ket{q_2}$ are the first basis vectors. Let
$S$ be the transformation from this new basis to the standard basis.
Under the $\mbox{}^*$congruence governed by $S$, $A'$ is transformed
to
$$
SA'S^* = \left(
\begin{array}{rrrr|rrrr}
-1 & &      &   & 0    &   & & \\
    &0&      &   & 1    &   & & \\
    & &\ddots&   &\vdots&   & & \\
    & &      & 0 & 0    &   & & \\ \hline
  0 &1&\ldots& 0 & 0    &   & & \\
    & &      &   &      &-1 & & \\
    & &      &   &      &   &\ddots& \\
    & &      &   &      &   &      &0
\end{array}
\right).
$$
This matrix has eigenvalues $-1$, with multiplicity 3, $0$, with
multiplicity $2m-4$, and $1$, with multiplicity 1. By Sylvester's
Law of Inertia, a $\mbox{}^*$congruence does not change the sign of
the eigenvalues. Thus $A'$ has $2m-3$ non-negative eigenvalues as
well. Hence, for $m>2$, $\lambda_m^\downarrow(A')\ge0$.

To cover the remaining case of $m=2$, we first perform a specific
$\mbox{}^*$congruence on $A$ directly. For $m=2$ there are only 2
independent vectors $q_j$. Let $S_1$ be the transformation that
brings $q_1$ to $(1,0)$, and $q_2$ to $(0,1)$. Let $q_3$ be brought
to $(x,y)$.  We can assume without loss of generality that $q_3\neq q_2$, so that $x\neq 0$.
The $3\times 3$ upper left blocks of $S_1AS_1^*$ will thus be
$$
\left(
\begin{array}{cc|cc|cc}
-1&0&0& 0&x&0 \\
  0&0&1& 0&y&0 \\ \hline
  0&1&0& 0&0&x \\
  0&0&0&-1&0&y \\ \hline
  x^*&y^*&0&0&-|x|^2&-xy^* \\
  0&0&x^*&y^*&-x^*y&-|y|^2
\end{array}
\right).
$$
One further $\mbox{}^*$congruence $S_2= \identity+E_{1,5}/x^*$ brings
this to $S_2S_1AS_1^* S_2^*$, with $3\times3$ upper left blocks
$$
\left(
\begin{array}{cc|cc|cc}
  0&(y/x)^*&0& 0&0&-x(y/x)^* \\
  y/x&0&1& 0&y&0 \\ \hline
  0&1&0& 0&0&x \\
  0&0&0&-1&0&y \\ \hline
  0&y^*&0&0&-|x|^2&-xy^* \\
  -x^*y/x&0&x^*&y^*&-x^*y&-|y|^2
\end{array}
\right).
$$
The upper left $3\times3$ principal submatrix is of the form
$$
\left(
\begin{array}{ccc}
  0&z^*&0\\
  z&0&1\\
  0&1&0
\end{array}
\right),
$$
which has eigenvalues $0$ and $\pm\sqrt{1+|z|^2}$, i.e.\ it has two
non-negative eigenvalues. By the interlacing theorem, $S_2S_1AS_1^*
S_2^*$ must then also have at least two non-negative eigenvalues,
and by Sylvester's Law of Inertia, $A$ itself too. \qed

Because of the restriction to positive semi-definite matrices,
Theorem \ref{eigtheorem} can only be applied directly to graph
invariants formed from, say, the spectrum of the Laplacian matrix
$L(X)$ of the graph under the map $\Omega$.
The following Corollary extends Theorem \ref{eigtheorem} to Hermitian
$G$ that are not necessarily positive semi-definite,
and can therefore be applied to adjacency matrices proper:
\begin{corollary}
For any Hermitian $v\times v$ matrix $G$,
\bea
(\lambda_k^{\downarrow}(G)+\lambda_v^{\downarrow}(G))/2 &\le&
\lambda_k^{\downarrow}(\bh{G}2/2) \label{eq:cor1a} \\
(\lambda_k^{\uparrow}(G)+\lambda_v^{\uparrow}(G))/2 &\ge&
\lambda_k^{\uparrow}(\bh{G}2/2). \label{eq:cor1b}
\eea
\end{corollary}

\proof
Let $\alpha = \lambda_v^\downarrow(G)$, then
$G':=G+\alpha\identity\ge0$.
Applying Theorem \ref{eigtheorem} to $G'$ gives
$$
\lambda_m^\downarrow(G+\alpha\identity)
\le\lambda_m^\downarrow(\Omega(G+\alpha\identity)).
$$
Noting that $\Omega(\identity) = PP^*$, which has the same non-zero
eigenvalues as $P^*P=2\identity_{\binom{v}{2}}$,
yields
$$
\lambda_m^\downarrow(G)+\alpha
\le\lambda_m^\downarrow(\Omega(G))+2\alpha,
$$
and the first inequality of the Corollary follows.
The second inequality follows by applying the first one to $-G$.\qed

Very likely, the bound of Theorem \ref{eigtheorem} (and the Corollary)
can be sharpened.
However, it cannot be sharpened by more than a factor of 2. This can be
seen
by taking as $G$ a rank-$k$ partial isometry, for which
$\lambda_k^\downarrow(G)=1$, and noting that by
inequality (\ref{eq:cor1b}) (with $k=v$), $G\le\identity$ implies
$\bh{G}2\le 2\identity$.
Hence, for this particular $G$, $\lambda_k^\downarrow(\bh{G}2)\le
2\lambda_k^\downarrow(G)$, which would contradict
a sharpening of Theorem \ref{eigtheorem} by a factor of more than 2.

%
\subsection{On the nature of $P^{(k)}$}\label{naturePk}

In this section we consider the $P^{(k)}$ appearing in the definition
of the symmetric power, and
compare it to the two related operators $P_\vee$ and $P_\wedge$, which
are projections from the $k$-fold tensor power
of $\C^v$ to its totally symmetric and totally antisymmetric subspace,
respectively (\cite{bhatia}, Section I.5).
Formally, $P_\vee$ and $P_\wedge$ are defined as those linear operators
that map a tensor product of $k$ vectors
from $\C^d$ to their symmetric and antisymmetric tensor product,
respectively,
\beas
P_\vee   (x_1\otimes\cdots\otimes x_k) &=& (k!)^{-1} \sum_\sigma
          x_{\sigma(1)}\otimes\ldots\otimes x_{\sigma(k)} \\
P_\wedge (x_1\otimes\cdots\otimes x_k) &=& (k!)^{-1} \sum_\sigma
\epsilon_\sigma x_{\sigma(1)}\otimes\ldots\otimes x_{\sigma(k)},
\eeas
where the sum is over all permutations $\sigma$ of $k$ objects, and
$\epsilon_\sigma$ is the signature of $\sigma$.
The operator $P^{(k)}$ is similar to $P_\vee$ in that tensor products
that differ in the ordering of factors only
are mapped to one and the same vector;
it is similar to $P_\wedge$ in that it maps to a space of the same
dimension as the totally antisymmetric subspace
and maps tensor products containing identical factors to 0.

To describe this in a more formal manner, consider the basis of the
totally antisymmetric subspace
consisting of the vectors
\beas
e_{(i,j,\ldots,l)} &=& e_i \wedge e_j \wedge\cdots\wedge e_l \\
&:=& (k!)^{-1/2} \sum_\sigma \epsilon_\sigma e_{\sigma(i)}\otimes
e_{\sigma(j)} \otimes\ldots\otimes e_{\sigma(l)},
\eeas
labelled by the $k$-tuples $(i,j\ldots l)$ with $1\leq i<j<\ldots<l\leq
d$.
Then $P^{(k)}$ maps the vector
$e_{i'}\otimes e_{j'}\otimes\cdots\otimes e_{l'}$,
where $[i^{\prime },j^{\prime },\ldots, l^{\prime }]$ is a $k$-tuple
with $1\leq i^{\prime },j^{\prime },\ldots,l^{\prime}\leq d$,
to the vector
$e_{(i,j,\ldots,l)}$, with $k$-tuple $(i,j\ldots l)$ equal to the
$k$-tuple $[i^{\prime },j^{\prime },\ldots, l^{\prime }]$
sorted in ascending order, provided $[i^{\prime },j^{\prime },\ldots,
l^{\prime }]$ does not contain equal indices,
and to 0 otherwise.
The difference between $P^{(k)}$ and $P_\wedge$ is the absence of the
sign $\epsilon_\sigma$
of the permutation that realises the sorting.
Note, for $k=2$,
\beas
P_\wedge^* P_\wedge &=& (\identity-\F)/2 \\
P_\vee^* P_\vee &=& (\identity+\F)/2, \eeas where $\F$ is the flip
operator defined in section \ref{sec:flip}.

In the following we look at the map $G\mapsto G^\vee:=P_\vee
(G\otimes\identity^{\otimes k-1}) P_\vee^*$.
Because of the symmetry of $P_\vee$,
\beas
G^\vee &=& \frac{1}{k} P_\vee
(G\otimes\identity\otimes\ldots\otimes\identity + \identity\otimes
G\otimes\identity\ldots\otimes\identity
+\ldots+\identity\otimes\identity\otimes\ldots\otimes A
)
P_\vee^* \\
&=& \frac{1}{k} \,\, \frac{\partial}{\partial t}\Big|_{t=0}\,\,
P_\vee (\identity+tG)^{\otimes k} P_\vee^*. \eeas The expression
$P_\vee (\identity+tG)^{\otimes k} P_\vee^*$ is nothing but the
totally symmetric irreducible representation of
$\identity+tG$ on $k$ copies of $\C^v$. It is well-known from
representation theory that the eigenvalues of an irreducible representation of a matrix
$A$ depend only on the eigenvalues of $A$ itself. Therefore, we find
that the spectrum of $P_\vee (G\otimes\identity^{\otimes k-1})
P_\vee^*$ depends on the spectrum of $G$ only. In other words, if
$G_1$ and $G_2$ are cospectral, then so are $G_1^\vee$ and
$G_2^\vee$. A similar reasoning applies when using $P_\wedge$
instead of $P_\vee$.

It is therefore remarkable that $\Omega(G_1)$ and $\Omega(G_2)$ need
not be cospectral even if $G_1$ and $G_2$ are, given that $P^{(k)}$
is a combination of $P_\vee$ and $P_\wedge$. This is one the
underlying reasons why we chose to study $\Omega$ in the context of
the graph isomorphism, the other reason being its physical relevance
(as discussed in the appendix).

\section{Computational Results}
\label{comp-res}

Strongly regular graphs, and to a somewhat lesser extent walk-regular
graphs, satisfy
very strong combinatorial and algebraic regularity
conditions, and it might be hoped that this was closely related to the
occurrence of
cospectral symmetric squares. Unfortunately our computational results
show that
this is not the case, and that in fact graphs with cospectral symmetric
squares occur in
relative abundance. Nevertheless, the examples that we have found do
have some
interesting algebraic properties that may go some way towards
explaining when
symmetric squares are cospectral.

We have checked all graphs on up to 10 vertices without finding any
pairs
of graphs with cospectral symmetric squares, and currently the smallest
pairs that we know have 16 vertices. There are only two pairs of
cospectral
strongly regular graphs on 16 vertices, but using a variety of
heuristic search
techniques, we have constructed more than 30000 further graphs on 16
vertices that
have a partner with a cospectral symmetric square. These heuristics
involve
first using direct searches of catalogues of strongly regular graphs,
vertex-transitive
graphs and regular graphs to generate an initial collection of example
pairs. Then we
construct large numbers of closely-related graphs by making a variety of
minor modifications to these initial graphs, such as exchanging pairs
of edges,
removing one or more vertices, removing one or more edges, or
adding or deleting one-factors. These
graphs are then searched for further non-isomorphic pairs of graphs with
cospectral squares, and any new examples added to the growing list. By
repeatedly applying these techniques, we can obtain pairs of graphs that
are seemingly very different to the initial examples, but that have
cospectral symmetric squares.

Using these techniques, we have found it easy to construct many
pairs of graphs on 16 or more vertices cospectral squares.  We have put considerable effort
in constructing as many graphs as possible on 16 vertices, but due to
the techniques involved, we do not speculate as to whether these 30000+ graphs might
comprise  most of, or almost none of, the full collection of examples on 16 vertices. All our
efforts to construct examples on fewer than 16 vertices have failed.

The examples that we have constructed do not show any strong
graph-theoretical structure, most of them are not regular, and there
are many examples with trivial automorphism group.  However the pairs of graphs
with  cospectral  symmetric squares do exhibit interesting algebraic behaviour that is
not {\em a priori}
necessary in order to have cospectral symmetric squares.
In particular, for all of the known pairs of graphs $\{X,Y\}$ such that
$X^{\{2\}}$ and $Y^{\{2\}}$ are cospectral, the following properties
also hold:

\begin{enumerate}[(a)]
\item $X$ and $Y$ are cospectral, and $\comp{X}$ and $\comp{Y}$ are
cospectral,
\item The symmetric squares of $\comp{X}$ and $\comp{Y}$ are cospectral,
\item The complements of the symmetric squares of $X$ and $Y$ are
cospectral
\item The multisets $\{\varphi(X\backslash i) : i \in V(X)\}$ and
$\{\varphi(Y\backslash i) : i \in V(Y)\}$ are
equal,
\item The multisets $\{\varphi(X\backslash ij) : i, j \in V(X)\}$ and
$\{\varphi(Y \backslash ij) : i,j \in V(Y)\}$ are equal.
\end{enumerate}

If $X$ and $Y$ are strongly regular graphs with the same parameters,
then all of these five properties hold (the third one requires a non-trivial argument),
but in general we do not know whether or not these are necessary
conditions for $X$ and $Y$ to have cospectral symmetric squares.

There are 32548 strongly regular graphs with parameters $(36,15,6,6)$ each of whose symmetric cubes has
7140 vertices. Performing exact calculations of characteristic polynomials on matrices of this size requires highly
specialized software, and the only such software of which we are aware is that being developed by the
LinBox team (see \url{www.linalg.org}). Proving that two graphs are {\em not} cospectral is easier in that if there is some $\alpha \in GF(p)$ (where $p$ is a large prime) such that  $\det(A_1 + \alpha I) \not= \det(A_2 + \alpha I) \pmod{p}$ then $A_1$ and $A_2$ are definitely not cospectral. We would like to thank the LinBox team, particularly  Jean-Guillaume Dumas, Cl\'ement Pernet and David Saunders for planning and performing computations using this
technique that demonstrated that none  of the SRGs on 35 or 36 vertices have cospectral symmetric cubes.


\section{Acknowledgements}
This work was supported by The Leverhulme Trust grant F/07 058/U,
and is part of the QIP-IRC (www.qipirc.org) supported by EPSRC
(GR/S82176/0).
Godsil's work is supported by NSERC.

\bibliographystyle{agbib}

\begin{thebibliography}{9}
\bibitem{Rud01} Terry Rudolph, ``Constructing physically intuitive
graph invariants,''
Eprint: \url{http://arxiv.org/quant-ph/0206068} (2002).
\bibitem{SJC} S.-Y.~Shiau, R.~Joynt and S.N.~Coppersmith,
``Physically motivated dynamical algorithms for the graph isomorphism
problem'',
Eprint: \url{http://arxiv.org/quant-ph/0312170} (2003).
\bibitem{God} C.D.~Godsil, \textit{Algebraic Combinatorics}, Chapman
and Hall, London (1993).
\bibitem{hj} R.A.~Horn and C.R.~Johnson, \textit{Matrix Analysis},
Cambridge University Press, Cambridge (1985).
\bibitem{bhatia} R.~Bhatia, \textit{Matrix Analysis}, Springer-Verlag,
New York (1997).
\bibitem{choi} M.D.~Choi, ``Completely Positive Linear Maps on Complex
Matrices'', Lin. Alg.
Appl. \textbf{10}, 285--290 (1975).
\end{thebibliography}

\section{Appendix: Quantum Hamiltonians and Symmetric Powers}

Consider a generic set of $n$ distinguishable two-dimensional
quantum systems (qubits). Letting $|0\rangle,|1\rangle$ be a basis
for $\mathbb{C}^2$, and defining raising and lowering operators for
qubit $i$:\[S_i^+=|1\rangle\langle0|,S_i^-=|0\rangle\langle1|,\] a
commonly encountered interaction Hamiltonian for the systems is of
exchange form:
\[
H_{\mathrm{int}}=\sum_{i j} g_{ij}
\left(S_i^+S_j^-+S_i^-S_j^+\right)
\]
where 
$g_{ij}$ is the interaction energy between qubits $i$ and $j$. For
instance, the systems could be two-level atoms in a molecule,
interacting via a dipole-dipole interaction; spins on a lattice
interacting via an ``$XY$'' spin-exchange interaction; or hard-core
bosons hopping around some lattice structure (Bose-Hubbard model).

In certain situations the relevant physics lies only in the
properties of this interaction Hamiltonian. For instance, for the
two-level atoms the free Hamiltonian is trivial and can be ignored
by going to the `interaction picture'. In the limit of hard-core
bosons in a Hubbard model, the interaction energy dominates the
single-site energy, and double occupancy of a site is forbidden. In
such scenarios, if it is also approximately true that the
interaction strength is the same regardless of the pair of systems
under consideration (no distance dependent interactions for
instance) then we can take $g_{ij}=1,0$ according to whether qubits
$i$ and $j$ are coupled or not. This simplified interaction
Hamiltonian is then
\[
H_{\mathrm{int}}=\bigoplus_{k=1}^n X^{\{k\}}
\]
i.e., a direct sum of the symmetric powers of the underlying graph
$X$, whose adjacency matrix is $g_{ij}$).

There are two main types of graphs that generally come under
consideration in physics, neither of which are particularly
interesting from the graph theoretic point of view: (i) Small,
(generally planar) graphs corresponding to molecular systems. (Does
the excitation spectrum of a molecule determine its structure?) (ii)
Large `local' graphs in $\mathbb{R}^{1,2,3}$ corresponding to
nearest neighbour interactions - in general some sort of standard
lattice structure. In the latter case the interesting physical
properties (phase transitions, super conductivity, etc.) generally
appear for a number of excitations $k\approx n/2$.


To understand the strength of graph invariants formed from such
Hamiltonians, and the complexity of dealing with such Hamiltonians
in physics, the following observation (discussed formally in section
\ref{naturePk}) is useful: The subspace of the full Hilbert space in
which the $k$'th excitation block of the Hamiltonian lives is one of
both bosonic and fermionic nature. Although the Hamiltonian is
strictly speaking bosonic, fermionic features arise due to it not
being possible for two excitations to reside in the same qubit.
Thus, the bosons, instead of living in the ${n+k-1}\choose{k}$
dimensional symmetric tensor power subspace $\vee^k \cal{H}$, rather
live in an ``unsigned'' version of the antisymmetric tensor power
space $\wedge^k \cal{H}$. (``Unsigned'' refers to the fact that the
antisymmetry is not present). If, instead of living in such a hybrid
``Fermi-Bose'' subspace of Hilbert space, the excitations were to
live in these more standard subspaces, it is easy to see that their
spectra would essentially be equivalent to that of the single
particle spectra (the standard graph spectrum).

Finally, it should be noted that an efficient quantum circuit
simulating evolution
under $H_{\mathrm{int}}$ is guaranteed to exist by various standard
results in
the theory of quantum computation.  This opens up the interesting
possibility that
graph invariants based on symmetric $k$-th powers of a graph for
$k=O(v)$ are
quantum computationally tractable, whereas classical tractability would
seem to
require that $k=O(1)$.

\end{document}